\baselineskip=15pt plus 2pt
\magnification =1200
\def\sqr#1#2{{\vcenter{\vbox{\hrule height.#2pt\hbox{\vrule width.#2pt
height#1pt\kern#1pt \vrule width.#2pt}\hrule height.#2pt}}}}
\def\square{\mathchoice\sqr64\sqr64\sqr{2.1}3\sqr{1.5}3}
 at 10truept
 at 10truept
 at 10truept
 at 10truept at 12truept

\font\medtenrm=cmr10 scaled\magstep2
\centerline {\medtenrm DISCRETE TRACY--WIDOM OPERATORS}\par
\vskip.1in
\centerline {GORDON BLOWER AND ANDREW MCCAFFERTY}\par
\vskip.05in
\centerline  {\sl Department of Mathematics and Statistics, 
Lancaster University}\par
\centerline  {\sl Lancaster, LA1 4YF, England, UK.}\par  
\centerline{ g.blower@lancaster.ac.uk, a.mccafferty@lancaster.ac.uk}\par
\vskip.05in
\centerline {(22nd November 2007)}\par
\centerline {( Revised 12th May 2008)}\par
\vskip.05in
\vskip.05in
\vskip.1in
{{\noindent {\sl Abstract} Integrable operators arise in random matrix theory,
where they describe the asymptotic eigenvalue distribution of large self-adjoint random
matrices from the generalized unitary ensembles. 
This paper considers discrete Tracy--Widom operators, and 
gives sufficient conditions
for a discrete integrable operator to be the square of a Hankel matrix.
Examples include the discrete Bessel kernel and kernels arising from
the almost Mathieu equation and the Fourier transform of Mathieu's
equation.}\par

\par

\noindent {\bf 1. Introduction}\par
\noindent We consider Tracy--Widom operators arising from 
first-order recurrence relations
$$ a(j+1)=T(j)a(j)\qquad (j=1, 2, \dots ),\eqno(1.1)$$
\noindent where $a(j)$ is a real $2\times 1$ vector and $T(x)$ is a
$2\times 2$ real matrix with entries that are rational functions of
$x$, and such that $\det T(j)=1$. Then with the skew-symmetric matrix
$$J=\left[\matrix{0&-1\cr 1&0\cr}\right]$$
\noindent the Tracy--Widom operator $K:\ell^2({\bf
N})\rightarrow \ell^2({\bf N})$ has matrix  
$$K(m,n)={{\langle Ja(m), a(n)\rangle }\over {m-n}}\qquad (m\neq
n)\eqno(1.2)$$
\noindent with respect to the usual orthonormal basis. In specific
examples, there are natural ways of defining the diagonal $K(n,n)$,
as we discuss below. We recall that a real sequence 
$(\phi_j)\in \ell^2({\bf N})$ gives a 
Hankel matrix $\Gamma_\phi =[\phi_{j+k-1}]_{j,k=1}^\infty;$
clearly $\Gamma_\phi$ is symmetric, and $\Gamma_\phi$ is Hilbert--Schmidt
if and only if $\sum_{n=1}^\infty n\vert \phi_n\vert^2$ converges. 
See [14].\par
\indent We consider whether such a $K$ may be expressed as
$K=\Gamma_\phi^2$, where $\Gamma_\phi$ is a Hankel matrix that is 
self-adjoint and Hilbert--Schmidt.\par 

\indent A significant example from random matrix theory is the discrete 
Bessel kernel
$$B(m,n)=\sqrt{\theta}{{{\hbox{J}}_m(2\sqrt{\theta})
{\hbox{J}}_{n-1}(2\sqrt{\theta})-{\hbox{J}}_{n}(2\sqrt{\theta}){\hbox{J}}_{m-1}(2\sqrt{\theta})}
\over{m-n}},\eqno(1.3)$$
\noindent as considered by Borodin {\sl et al} [5] and Johansson [11]. They showed that
$B$ is the square of the Hilbert--Schmidt Hankel matrix
$[{\hbox{J}}_{m+k-1}(2\sqrt{\theta })]$, and thus obtained information about
the spectrum of $B$ itself.\par

\indent Tracy and Widom observed that many of the fundamentally important
kernels in random matrix theory have the form
$$W(x,y)={{f(x)g(y)-f(y)g(x)}\over {x-y}}\qquad (x\neq y)\eqno(1.4)$$
\noindent where $f,g$ are bounded real functions in 
$L^2(0, \infty )$ such that 
$${{d}\over{dx}}\left[\matrix {f(x)\cr g(x)\cr}\right] =
\left[\matrix {\alpha (x)&\beta (x)\cr -\gamma (x)&-\alpha
(x)\cr}\right]\left[\matrix {f(x)\cr g(x)\cr}\right],\eqno(1.5)$$
\noindent and $\alpha (x), \beta (x)$ and $\gamma (x)$ are real
rational functions; see [16]. Then there is a bounded linear operator
$W:L^2(0, \infty )\rightarrow L^2(0, \infty )$ given by
$Wh(x)=\int_0^\infty W(x,y)h(y)dy$. Of particular importance is the 
case in which $W$ is a trace-class
operator such that $0\leq W\leq I$, since such a $W$ is associated with a
determinantal point process as in [15]. In order to verify this property in
special cases, Tracy and Widom showed that $W$ is the square of a
Hankel operator $\Gamma$ which is self-adjoint and Hilbert--Schmidt.
The spectral theory and realization of Hankel operators is well 
understood [8, 13], so such a
factorization is valuable. In [3, 4], we developed this method further 
and considered which differential
equations lead to kernels that can be factored as squares of Hankel
operators. The factorization theorems involve the notion of operator
monotone functions.\par

\indent Here we consider various discrete Tracy--Widom kernels
and their factorization as Hankel products, using the formal 
analogy between differential equations and 
difference equations which suggests likely factorization theorems. We
can write the matrix in (1.5) as $J\Omega (x)$, where $\Omega (x)$ is
real and symmetric, and then consider the analogous one-step
transition matrix to be $T(x)=\exp (J\Omega (x)).$ The
functions $x\mapsto x$ and $\mapsto -1/x$ are operator monotone monotone
increasing on $(0, \infty )$, and they appear in the transition
matrices for the discrete Bessel kernel in section 3 and the discrete analogue of
the Laguerre kernel in section 4.\par

\indent The functions $x\mapsto x^2$ and $x\mapsto -1/x^2$ are not
operator monotone increasing by [9], and so we cannot hope to have simple
factorization theorems when they appear in the transition matrix $T$.
In the special case of the parabolic cylinder equation
$-\phi_n''(x)+(x^2/4-1/2)\phi_n (x)=n\phi_n(x)$, Aubrun [1] recovered a
factorization of the corresponding kernel
$$K(x,y)={{\phi_n(x)\phi_{n-1}(y)-\phi_{n-1}(x)\phi_n(y)}\over {x-y}}$$
\noindent in the form 
$$K=\Gamma_\phi\Gamma_\psi +\Gamma_\psi \Gamma_\phi\eqno(1.6)$$
\noindent where $\Gamma_\phi$ and $\Gamma_\psi$ are bounded and
self-adjoint Hankel operators. The parabolic cylinder function is a
confluent form of Mathieu's functions associated with the elliptic
cylinder, since the parabola is the
limiting case of an ellipse as the eccentricity increases to one [17,
p. 427]. Hence
it is natural to factorize kernels associated with Mathieu's equation 
$${{d^2u}\over{dz^2}}+(\alpha +\beta \cos z) u(z)=0\eqno(1.7)$$
\noindent in the form (1.6). In section 5 we consider the first-order difference equation associated with the Fourier
transform of Mathieu's equation. In section 6, we consider almost
Mathieu operators.\par
\indent For a compact and self-adjoint operator $W$, the spectrum
consists of real eigenvalues $\lambda_j$ which may be ordered so
that the sequence of singular numbers $s_j=\vert \lambda_j\vert$
satisfies $s_1\geq s_2\geq \dots $. While the factorization $W=\Gamma^2$ immediately determines the
spectrum of $W$ from the spectrum of $\Gamma$, a factorization (1.6)
imposes bounds upon the singular numbers of $K$ is terms of the
eigenvalues of $\Gamma_\phi$ and $\Gamma_\psi$.  
\vskip.1in

\noindent {\bf 2. Factoring discrete Tracy--Widom operators as squares of Hankel
matrices}\par

\indent {\bf Theorem 2.1.} {\sl Let $T(j)$ and $B(j)$ be $2\times 2$
real matrices, and $(a(j))$ a sequence of real $2\times 1$ vectors such
that 
$$a(j+1)=T(j)a(j)\qquad (j\in{\bf N}), \eqno(2.1)$$
$$a(j)\rightarrow 0\qquad (j\rightarrow\infty),\eqno(2.2)$$
$$\sum_{j=1}^\infty\Vert B(j)a(j)\Vert^2 <\infty .\eqno(2.3)$$
\noindent Suppose further that there exists a real symmetric matrix $C$
with eigenvalues $0$ and $\lambda$, where $\lambda <0$, such that 
$${{T(n)^tJT(m)-J}\over{m-n}}=B(n)^tCB(m)\qquad (m\neq n;
m,n \in {\bf N}).\eqno(2.4)$$
\indent Let $\phi(j)=\vert\lambda\vert^{1/2}\langle v_\lambda , B(j)a(j)\rangle$ where 
$v_\lambda$ is a real unit eigenvector corresponding to
$\lambda$. Then  $\Gamma_\phi $ is compact and 
$K=\Gamma_\phi^2$ has entries}
$$K(m,n)={{\langle Ja(m), a(n)\rangle }\over {m-n}}\qquad (m\neq
n; m,n\in{\bf N}).\eqno(2.5)$$
\vskip.05in
\indent {\bf Proof.} Let $K(m,n)$ be as in (2.5) and let 
$$G(m,n)=K(m,n)
-\sum_{k=1}^\infty\phi(m+k-1)\phi(n+k-1),\qquad (m\neq
n;m,n\in{\bf N})\eqno(2.6)$$ where the
infinite sum converges because of the condition (2.3). Observe that $a(j)\rightarrow 0$ implies that
the first term in $G(m,n)$ tends to zero as $m$ or $n\rightarrow\infty$, and that the same is
true of the Hankel sum:
$$
\Bigl\vert
\sum_{k=1}^{\infty}\phi(m+k-1)\phi(n+k-1)\Bigr\vert\leq\left(\sum_{k=1}^\infty\vert\phi(m+k-1)\vert^2\right)^{1/2}\left(\sum_{k=1}^\infty\vert\phi(n+k-1)\vert^2\right)^{1/2}\rightarrow
0,\eqno(2.7)
$$
as $m$ or $n\rightarrow\infty$.
Hence $G(m,n)\rightarrow 0$ as $m$ or $n\rightarrow\infty$. Now let $U$ be the real orthogonal matrix with $v_\lambda$ in the first
column, and the eigenvector of $C$ corresponding to the eigenvalue $0$ in
the second column. Then $U^tCU={\hbox{diag}}(\lambda,0)$, and we have
$$\eqalignno{K(m+1,n+1)-K(m,n)&={{1}\over {m-n}}\bigg(\bigl\langle
JT(m)a(m),T(n)a(n)\bigr\rangle-\bigl\langle
Ja(m),a(n)\bigr\rangle\bigg)\cr
&={{1}\over{m-n}}\left\langle(T(n)^tJT(m)-J)a(m),a(n)\right\rangle\cr
&=\left\langle B(n)^tCB(m)a(m),a(n)\right\rangle\cr
&=\left\langle U {\hbox{diag}}(\lambda,0)U^tB(m)a(m),B(n)a(n)\right\rangle\cr
&=\lambda\left\langle{\hbox{diag}}(1,0)U^tB(m)a(m),{\hbox{diag}}(1,0)U^tB(n)a(n)\right\rangle\cr
&=-\phi(m)\phi(n).&(2.8)}
$$
The above calculation, and the equality
$$
\sum_{k=1}^{\infty}\phi(m+k)\phi(n+k)-\phi(m+k-1)\phi(n+k-1)=-\phi(m)\phi(n)
\eqno(2.9)$$
together imply that
$G(m+1,n+1)=G(m,n)$, and so in fact $G(m,n)=0$ for all $m,
n\in{\bf N}$,
which gives 
$$K(m,n)=\sum_{k=1}^\infty \phi (m+k-1)\phi (n+k-1)\qquad (m\neq n;m,
n\in{\bf N}),\eqno(2.10)$$
\noindent where the right-hand side is the $(m,n)^{th}$ entry of the square of a Hankel
matrix $\Gamma_\phi.$ Further, we observe that $K$ is the composition of the Hilbert transform with
the compact operators\par
\noindent ${\ell^2}({\bf N}; {\bf
C})\rightarrow {\ell^2}({\bf N}; {\bf C}^2)$ given by $(x_n)\mapsto (Ja(n)x_n)$ and the adjoint of $(x_n)\mapsto (a(n)x_n)$; 
so $K$ is compact. We deduce that $\Gamma_\phi$ is also compact.\par
\rightline{$\square$}\par
\vskip.1in 
\indent {\bf Proposition 2.2.} {\sl 
Let $a(j)$, $T(j)$ and $B(j)$
satisfy conditions (2.1) and  (2.2), and suppose further that 
$$\sum_{j=1}^\infty j\Vert B(j)a(j)\Vert^2 <\infty .\eqno(2.11)$$
\noindent Now let $K=\Gamma_\phi^2$ be as in Theorem 2.1.\par
\indent (i) Then $K$ is a positive semidefinite and trace-class 
operator.}\par
\indent {\sl (ii) For each $n\in {\bf N},$ there exist self-adjoint Hankel operators
 $\Gamma_n$, where $\Gamma_n$ has rank at most $n$, such that}
$$s_n(K)  =\Vert \Gamma_\phi -\Gamma_n\Vert^2,\eqno(2.12)$$
\noindent {\sl so $\Gamma_n^2\rightarrow K$ as
$n\rightarrow\infty$.}\par
\vskip.05in
\indent {\bf Proof.} (i) The Hilbert--Schmidt norm of $\Gamma_\phi$ satisfies
$\Vert\Gamma_\phi\Vert^2_{HS}=\sum_{k=1}^\infty k\phi (k)^2<\infty
.$ Hence $K=\Gamma_\phi^2$ is of trace class.\par
\indent (ii) Since $\Gamma_\phi$ is self-adjoint, the singular numbers
satisfy
$s_n(K)=s_n(\Gamma_\phi^2 )=s_n(\Gamma_\phi)^2$. By the Adamyan--Arov--Krein
theorem [14], there exists a unique Hankel operator $\Gamma_n$ with rank at
most $n$ such that $s_n(\Gamma_\phi )=\Vert \Gamma_\phi -\Gamma_n\Vert .$
Evidently $\Gamma_n^*$ is also a Hankel operator of rank at most $n$
such that $s_n(\Gamma_\phi )=\Vert \Gamma_\phi -\Gamma^*_n\Vert$, so by
uniqueness $\Gamma_n=\Gamma_n^*$.\par
\indent We have $\Vert\Gamma_n-\Gamma_\phi\Vert\rightarrow 0$ as
$n\rightarrow\infty$, so $\Gamma_n^2\rightarrow \Gamma_\phi^2$ as
$n\rightarrow\infty.$\par 
\rightline{$\square$}\par

\indent {\bf Definition.} For a compact 
and self-adjoint operator $W$ on a Hilbert space $H$, the spectral
multiplicity function $\nu_W:{\bf R}\rightarrow
\{0,1,\dots\}\cup\{\infty\}$ is given by
$$
\nu_W(\lambda)={\hbox{dim}}\{x\in H:Wx=\lambda x\}\quad 
(\lambda\in{\bf R}).
\eqno(2.13)$$
\par
\indent {\bf Proposition 2.3. }{\sl Let $K$ be as in Proposition 2.2.
Then the following hold:\par
\indent (i) $\nu_{K}(0)=0 {\hbox{ or }}\nu_K(0)=\infty;$\par
\indent (ii) $\nu_K(\lambda)<\infty$ and
  $\nu_K(\lambda)=\nu_{\Gamma_\phi}(\sqrt{\lambda})
+\nu_{\Gamma_\phi}(-\sqrt{\lambda})$ for all $\lambda >0$.\par
\indent (iii) If $\nu_K(\lambda)$ is even, then 
 $\nu_{\Gamma_\phi}(\sqrt{\lambda})=\nu_{\Gamma_\phi}
(-\sqrt{\lambda}).$\par
\indent (iv) If $\nu_K(\lambda)$ is odd, then 
 $\vert\nu_{\Gamma_\phi}(\sqrt{\lambda})-\nu_{\Gamma_\phi}
(-\sqrt{\lambda})\vert=1.$}
\vskip.05in
\indent {\bf Proof. }
{\sl (i)} follows from Beurling's theorem (see [14, p.15]), while {\sl (ii)} is
elementary. Peller, Megretski\u \i\quad and Treil show in [13] that for
any compact and self-adjoint Hankel operator $\Gamma_\phi$, the spectral
multiplicity function satisfies
$\vert\nu_{\Gamma_\phi}(\lambda)-\nu_{\Gamma_\phi}
(-\lambda)\vert\leq 1$. Using this, and
{\sl (ii)}, statements
{\sl (iii)} and {\sl (iv)} follow immediately.
\vskip.05in

\noindent {\bf 3. The Discrete Bessel Kernel}\par
\noindent We show how Theorem 2.1 can be applied to the discrete Bessel kernel to
recover a result from [5, 11].\par
\vskip.05 in
\indent {\bf Proposition 3.1.}
{\sl Let} ${\hbox{J}}_n(z)$ {\sl be the Bessel function of the first 
kind of order
$n$, let} ${\hbox{J}}_n={\hbox{J}}_n(2\sqrt{\theta})$, {\sl where
$\theta>0$; let} $\phi(n)={\hbox{J}}_{n+1}$ {\sl and}
$a(n)=[\sqrt{\theta}{\hbox{J}}_n,{\hbox{J}}_{n+1}]^t$. {\sl Then the Hankel operator 
$\Gamma_\phi$
is Hilbert--Schmidt, and $B=\Gamma_\phi^2$ has entries}
$$B(m,n)={{\langle Ja(m),a(n)\rangle}\over{m-n}} 
\quad(m\neq n; m,n\in{\bf N}).\eqno(3.1)$$}
\par
\indent {\bf Proof. }
It is clear that (2.1) holds, since we have the recurrence relation
$$
{\hbox{J}}_{n+2}(2z)={{n+1}\over{z}}{\hbox{J}}_{n+1}(2z)-{\hbox{J}}_{n}(2z),
\eqno(3.2)$$
giving $a(n+1)=T(n)a(n)$, where
$$
T(n)=\left[\matrix{ 0&\sqrt{\theta}\cr
{{-1}\over{\sqrt{\theta}}}& {{n+1}\over{\sqrt{\theta}}}\cr}\right].
\eqno(3.3)$$
Note that 
$$
{{T(n)^tJT(m)-J}\over{m-n}}=C,
\eqno(3.4)$$
where $C={\hbox{diag}}(0,-1)$, which is clearly of rank one. The non-zero eigenvalue of $C$ is
$\lambda=-1$, and a corresponding unit eigenvector is
$v_\lambda=[0,1]^t$, so 
$$\vert\lambda\vert^{1/2}\langle
v_\lambda,a(n)\rangle={\hbox{J}}_{n+1}=\phi(n).$$ We now verify condition
(2.11), and thus (2.2). Note that
$${{1}\over{\theta}}
\sum_{n=1}^\infty
n{\hbox{J}}_{n+1}^2<{{1}\over{\theta}}\sum_{n=1}^\infty
(n+1)^2{\hbox{J}}_{n+1}^2=\sum_{n=1}^\infty({\hbox{J}}_{n+2}
+{\hbox{J}}_n)^2\leq 4\sum_{n=1}^\infty {\hbox{J}}_n^2.
\eqno(3.5)$$
The
standard formula (see [17, p.379])
$$e^{i2\sqrt{\theta}\sin{\psi}}
={\hbox{J}}_0(2\sqrt{\theta})+2\sum_{m=1}^{\infty}{\hbox{J}}_{2m}(2\sqrt{\theta})
\cos 2m\psi+
2i\sum_{m=1}^{\infty}{\hbox{J}}_{2m-1}(2\sqrt{\theta} )\sin(2m-1)\psi
\eqno(3.6)$$
and Parseval's identity can be used to show that
${\hbox{J}}_0(2\sqrt{\theta})^2+2\sum_{m=1}^{\infty}{\hbox{J}}_{m}
(2\sqrt{\theta})^2=1$ for all $\theta >0$, and
hence that the sum on the right hand side of (3.5) is finite.

\vskip.1in  
\noindent {\bf 4. A discrete analogue of the Laguerre differential
equation}\par
\vskip.05in
\noindent In this section we consider a case in which condition (2.2) is
violated and we cannot hope to factor the kernel $K$ as the square of a
Hankel operator. Nevertheless, we can identify a Toeplitz operator $W$
such that $K-W$ factors as a product of Hankels.\par
\vskip.05in

\indent {\bf Proposition 4.1.} {\sl For $\theta\in {\bf R}$, let $(a(j))$ satisfy the recurrence relation $a(j+1)=T(j)a(j)$
with 
$$T(j)=\left[\matrix{ \theta/(j+1)& -1\cr 1& 0\cr}\right]\eqno(4.1)$$
\noindent and $a(1)=[\theta ,1]^t$. Then there exist polynomials
$p_j(\theta)$ of degree $j$ such that\par
\indent (i)  $a(j)=[p_j(\theta ),
p_{j-1}(\theta )]^t$.\par 
\indent (ii) The self-adjoint Hankel matrix $\Gamma_\phi =[\phi
(j+k-1)]_{j,k=1}^\infty$ with entries
$$\phi (j)={{p_j(\theta )}\over {j+1}}\eqno(4.2)$$
\noindent is a bounded linear operator such that 
$\theta \Gamma_\phi^2=K+W$ where $K$ has entries}
$$K(m,n)={{\langle Ja(m), a(n)\rangle}\over{m-n}}\qquad (m\neq n;
m,n\in{\bf N} ),\eqno(4.3)$$
\noindent {\sl and $W$ is a bounded Toeplitz operator with matrix
$$W(m-n)={{\langle J^{\bar m-\bar n +1}T_\infty a(1), T_\infty
a(1)\rangle}\over {m-n}}\qquad(m\neq n)$$
\noindent for some $2\times 2$ matrix $T_\infty$, and $\bar m$ and $\bar n$ are the congruence classes of $m$ and
$n$ modulo $4$.}\par

\vskip.05in
\indent {\bf Proof.} When $\theta=0$, the recurrence relation reduces
to $a(n+1)=Ja(n)$, with solution $a(n)=J^{n-1}a(1)$. This gives rise to
a kernel 
$$K(m,n)={{\langle Ja(m), a(n)\rangle}\over{m-n}}={{\langle J^{m-n+1}a(1),
a(1)\rangle}\over {m-n}},\qquad (m\neq n)\eqno(4.4)$$
\noindent which has the shape of a Toeplitz operator, and is a variant on the Hilbert transform matrix
$[1/(m-n)]_{m\neq n}.$\par
\indent Now suppose $\theta\neq 0$. The matrix $J$ satisfies $J^4=I$, and so we consider the
partial product of the $T(j)$ in bunches of four, with the $j^{th}$
bunch giving
$$B(j)=T(4j)T(4j-1)T(4j-2)T(4j-3)=I-{{\theta}\over{2j}}J+O(1/j^2)\quad
(j\in{\bf N}).\eqno
(4.5)$$  
\noindent Now we deduce that 
$$\Vert B(j)\Vert^2=\Vert B(j)^*B(j)\Vert =1+O(1/j^2),\eqno(4.6)$$
\noindent and likewise with $B(j)^{-1}$ in place of $B(j)$; so there
exists $C(\theta )$ such that 
$$\Vert T(n)T(n-1)\dots T(2)T(1)\Vert\leq C(\theta ),$$
$$\Vert T(1)^{-1}T(2)^{-1}\dots T(n-1)^{-1}T(n)^{-1}\Vert\leq
C(\theta )\qquad (n=1,2, \dots ).\eqno(4.7)$$
\noindent It follows that there exists $\kappa (\theta )>0$ such
that $\kappa (\theta )<\Vert a(n)\Vert<\kappa (\theta )^{-1}$ for all
$n$, so (2.2) is violated. We introduce
$$C_k=\exp\Bigl(\theta J\sum_{j=1}^k{{1}\over{2j}}\Bigr)B(k)B(k-1)\dots
B(1),\eqno(4.8)$$
\noindent which satisfies $C_{k+1}-C_k=O(1/k^2)$; so the limit
$$T_\infty =\lim_{k\rightarrow\infty }C_k\eqno(4.9)$$
\noindent exists. One can check that 
$$\exp\Bigl(\sum_{j=1}^{k+n}{{\theta J}\over{2j}}\Bigr)^*
\exp\Bigl(\sum_{j=1}^{k+m}{{\theta J}\over{2j}}\Bigr)\rightarrow
I\eqno(4.10)$$
\noindent as $k\rightarrow\infty$, and hence
$$\eqalignno{\langle Ja(m+4k),& a(n+4k)\rangle\cr
& =\langle JT(m+4k)T(m+4k-1)\dots
T(1)a(1), T(n+4k)T(n+4k-1)\dots T(1)a(1)\rangle\cr
&\rightarrow \langle J^{\bar m+1}T_\infty a(1), 
J^{\bar n}T_\infty a(1)\rangle\qquad (k\rightarrow\infty ).&(4.11)\cr}$$
\indent  For temporary convenience we introduce
$$\tilde K(m,n)=\cases{ {{\langle Ja(m), a(n)\rangle}\over{
(m-n)}}& for
$m\neq n$; \cr
0& for $m=n$.\cr}\eqno(4.12)$$
\noindent The discrete Hilbert transform is bounded on $\ell^2({\bf
Z}; {\bf C}^2)$ by [7], so $\tilde K$ defines a
bounded linear operator on $\ell^2({\bf N})$, but condition
(2.2) is violated.\par 
\indent The $p_j(\theta )$ satisfy the recurrence
relation 
$$p_{n+1}(\theta )+p_{n-1}(\theta )={{\theta}\over{n+1}}p_n(\theta
)\eqno(4.13)$$
\noindent with $p_0(\theta )=1$ and $p_1(\theta )=\theta$, so clearly
$p_j(\theta )$ is a polynomial of degree $j$ such that $\vert
p_n(\theta )\vert\leq \kappa (\theta )^{-1}$ for all $\theta$ and $n$. Further, 
$${{T(n)^tJT(m)-J}\over {m-n}}=\left[\matrix{-\theta/(m+1)(n+1)&0\cr 
0& 0\cr}\right];\eqno(4.14)$$
\noindent so that, by the calculation in the proof of Theorem 2.1,
$$\tilde K(m+1, n+1)-\tilde K(m,n)=
-\theta {{p_m(\theta )p_n(\theta )}\over {(m+1)(n+1)}}.\eqno(4.15)$$
\indent We can write 
$$\tilde K(m,n)-\tilde K(m+N+1, n+N+1)=\theta \sum_{k=0}^N {{p_{m+k}(\theta )p_{n+k}(\theta
)}\over {(m+k+1)(n+k+1)}}\qquad (m\neq n)\eqno(4.16)$$
\noindent where the limit 
$$\lim_{N\rightarrow\infty}  \tilde K(m+N+1, n+N+1)=W(m-n)\qquad
(m\neq n)\eqno(4.17)$$
\noindent exists and is finite since the sequence $\phi (j)
=p_j(\theta )/(j+1)$ is square summable; so 
$$\tilde K(m,n)=\theta \sum_{k=0}^\infty {{p_{m+k}(\theta )p_{n+k}(\theta )}\over
{(m+k+1)(n+k+1)}} +W(m-n)\qquad (m\neq n).\eqno(4.18)$$
\noindent We now define $W(0)=0$, and let 
$$ K(m,n)=\theta \sum_{k=0}^\infty {{p_{m+k}(\theta )p_{n+k}(\theta )}\over
{(m+k+1)(n+k+1)}} +W(m-n);\eqno(4.19)$$
\noindent so that the matrix of $K$ equals the matrix of $\tilde K$,
except on the principal diagonal, and the principal diagonal of $K$ is
a bounded sequence; hence $K$ is a bounded
linear operator and satisfies the preceding
identities also for $m=n$. Let $S$ be the shift operator on
$\ell^2({\bf N})$. Now $\theta \Gamma_{\phi}^2$ equals the 
limit in the weak operator topology
of the sequence $K-S^{*n}KS^{n}$ as $n\rightarrow\infty$, so $\Gamma_{\phi}$ is
bounded and hence $W=K-\theta\Gamma_{\phi}^2$ is also bounded. We
recognise the matrix of $W$ as 
$$W(m-n)={{\langle J^{\bar m-\bar n +1}T_\infty a(1), T_\infty
a(1)\rangle}\over {m-n}}\qquad (m\neq n).$$
\rightline{$\square$}\par 
\vskip.05in
\indent {\bf Remark 4.2.} The generating function $f(z)=\sum_{j=0}^\infty p_j(\theta
)z^j$ satisfies the differential equation
$$(1+z^2)f'(z) +(2z-\theta )f(z)=0$$
\noindent with initial condition $f(0)=1$, and hence 
$$f(z)=\Bigl({{1-iz}\over{1+iz}}\Bigr)^{i\theta /2}{{1}\over
{1+z^2}}.\eqno(4.20)$$
\indent For comparison, Laguerre's equation [16] may be expressed as 
$${{d}\over{dx}}\left[\matrix{u(x)\cr u'(x)\cr}\right]=
\left[\matrix{0&1\cr {{1}/{4}}-{{(n+1)}/{x}}&0\cr}\right]
\left[\matrix{u(x)\cr u'(x)\cr}\right],
\eqno(4.21)$$ 
\noindent with solution $u(x)=xe^{-x/2}L_n^{(1)}(x)$ where 
$$L_n^{(1)}(x)={{x^{-1}e^x}\over {n!}} {{d^n}\over{dx^n}}\bigl(
x^{n+1}e^{-x}\bigr)\qquad (x>0)
\eqno(4.22)$$
\noindent is the Laguerre polynomial of degree $n$ and parameter 
$1$. The Laplace transform of $u$ is the rational function
$${\cal L}(u; \lambda )=(n+1){{(\lambda -{{1}\over{2}})^n}
\over {(\lambda +{{1}\over{2}})^{n+2}}}\qquad
(\Re\lambda >-1/2).\eqno(4.23)$$

\vskip.05in

\noindent {\bf 5. The Fourier transform of Mathieu's equation}\par
\noindent Let $\beta\neq 0$ be a real number; then there exists
a sequence of real values of $\alpha$ such that Mathieu's equation
$${{d^2u}\over{d\theta^2}}+(\alpha +\beta \cos\theta )u(\theta
)=0\eqno(5.1)$$
\noindent has a real periodic solution with period $2\pi$ or
$4\pi$. The odd or even periodic solutions are known as Mathieu
functions, and various determinants describe the dependence of the 
eigenvalues $\alpha$ on $\beta$, as in [12, 17]. Here we are
concerned with some matrices that arise from the Fourier transform of
the differential equation.\par
\vskip.05in
\indent {\bf Theorem 5.1.} {\sl Suppose that $u$ has Fourier expansion 
$u(\theta )=\sum_{n=-\infty}^\infty b_ne^{in\theta }.$ 
Let $\Gamma_u$ be the Hankel matrix $[b_{j+k-1}]_{j,k=1}^\infty$, 
let $\Gamma_v$ be the Hankel matrix $[(j+k-1)b_{j+k-1}]_{j,k=1}^\infty$ and let 
$K=(-2/\beta )(\Gamma_u\Gamma_v+\Gamma_v\Gamma_u).$ Then $K$ is a trace
 class operator on $\ell^2({\bf N})$ such that}
$${\hbox{trace}}(K)={{1}\over{\beta\pi }}\int_0^{2\pi}\bigl\vert
{{du}\over {d\theta}}\bigr\vert^2\, d\theta $$ 
\noindent {\sl and}
$$K(j,k)={{b_{j-1}b_{k}-b_{j}b_{k-1}}\over {j-k}}\qquad 
(j,k\in{\bf N}; j\neq k).\eqno(5.2)$$
\vskip.05in
\indent {\bf Proof.} Since $u$ is real, we have
$b_m=\overline{b_{-m}}$. The recurrence relation for the Fourier 
coefficients 
$$2(-n^2+\alpha )b_n+\beta b_{n+1}+\beta b_{n-1}=0\eqno(5.3)$$
may be expressed as the first-order recurrence relation
$$\left[\matrix{b_{n}\cr b_{n+1}\cr}\right]=
\left[\matrix{  0&1\cr -1& (2/\beta )(n^2-\alpha )\cr}\right]
\left[\matrix{b_{n-1}\cr b_{n}\cr}\right],\eqno(5.4)$$
\noindent or in the obvious shorthand $a(n+1)=T(n)a(n).$ Then we have 
$${{T(n)^tJT(m)-J}\over{m-n}}=\left[\matrix{ 0&0\cr
0&(-2/\beta)(m+n)\cr}\right].\eqno(5.5)$$
\noindent We introduce the kernel $\tilde K$ by the formula
$$\tilde K(m,n)={{\langle Ja(m),
a(n)\rangle}\over{m-n}}={{b_{n}b_{m-1}-b_{n-1}b_{m}}\over {m-n}}$$
\noindent which therefore satisfies 
$$\tilde K(m+1, n+1)-\tilde K(m,n)=(-2/\beta
)(mb_mb_n+nb_nb_m)\eqno(5.6)$$
\noindent and $\tilde K(m,n)\rightarrow 0$ as $m,n\rightarrow\infty$ in any way such that $m\neq n.$ We deduce that
$$\tilde K(m,n)={{2}\over{\beta }}\sum_{k=0}^\infty
(m+k)b_{m+k}b_{n+k}+(n+k)b_{n+k}b_{m+k}\qquad (m,n\in{\bf N}; m\neq n).$$ 
\noindent Hence $\tilde K(m,n)$ is the $(m,n)$ entry of the matrix of 
$K=(2/\beta)(\Gamma_u\Gamma_v+\Gamma_v\Gamma_u)$ for all $m\neq n.$ \par
\indent Since $u$ and $u''$ are square integrable, 
the series $\sum_{n=-\infty}^\infty n^4\vert b_n\vert^2$ converges; 
so $\Gamma_u$ and $\Gamma_v$ are Hilbert--Schmidt, and $K$ is trace 
class. Further, we have
$$\eqalignno{{\hbox{trace}}(K)&=\sum_{m=1}^\infty K(m,m)\cr
&={{4}\over{\beta}}\sum_{m,k=1}^\infty (m+k-1)b_{m+k-1}^2\cr
&={{4}\over{\beta}}\sum_{m=1}^\infty m^2b_m^2\cr
&={{2}\over{\beta}}\int_0^{2\pi} \vert u'(\theta )\vert^2 \,
{{d\theta}\over{2\pi}} .
&(5.7)\cr}$$   
\vskip.1in

\noindent {\bf 6. Almost Mathieu operators}\par
\noindent We introduce the almost Mathieu operator $H:\ell^2({\bf
Z})\rightarrow \ell^2({\bf Z})$ by
$$(Hu)_n=u_{n+1}+u_{n-1}+\lambda \cos 2\pi(n\theta
+\alpha )\,u_n\eqno(6.1)$$
\noindent for $u=(u_n)_{n=-\infty}^\infty\in \ell^2({\bf
Z})$, where $(Hu)_n$ denotes the $n^{th}$ term in the sequence $Hu\in \ell^2({\bf
Z})$. For all real $\lambda ,\theta, \omega$ and $\alpha$, the
operator $H$ is bounded and self-adjoint, with spectrum contained in
$[-2-\vert\lambda\vert , 2+\vert\lambda\vert ]$. According to the
precise values of the parameters, as we discuss below, the spectrum can consist of a mixture
of point spectrum, continuous spectrum and singular continuous
spectrum.\par  
\indent {\bf Definition.} Let 
$E$ be an eigenvalue of $H$ with the corresponding
eigenvector $(u_n)$. Say that $(u_n)$ has {\sl exponential decay} if 
there exist $C,\delta >0$ and $n_0\in {\bf Z}$ such that
$$\vert u_n\vert\leq Ce^{-\delta \vert n-n_0\vert}\qquad (n\in {\bf
Z}).\eqno(6.2)$$
\noindent (Typically, $\delta$ depends upon $\theta, \alpha$ and
$\lambda$ in a complicated fashion.) Say that $H$ exhibits {\sl Anderson
localization} if its spectrum is pure point and all eigenvectors are of
exponential decay. \par  
\vskip.1in

\indent {\bf Definition.} Say that $\theta\in {\bf R}$ is
{\sl Diophantine}
if there exist $c(\theta )>0$ and $r(\theta )>0$ such that 
$$\vert \sin 2\pi j\theta\vert\geq c(\theta )\vert j\vert^{-r(\theta
)}\qquad (j\in {\bf Z}\setminus \{ 0\} ).\eqno(6.3)$$
\indent With respect to Lebesgue measure, almost all real numbers are
Diophantine; see [2, p. 373]. Clearly rational numbers are not Diophantine, nor are
Liouville numbers; see [10].\par

\indent Following the work of several mathematicians, as summarized in 
[6], Jitomirskaya [10]
obtained a satisfactory description of the spectrum of the almost
Mathieu operator.\par
\vskip.05in
\indent {\bf Lemma 6.1.} [10] {\sl Suppose that $\theta$ is Diophantine.
Then there exists a set $S_\theta$ such that ${\bf R}\setminus S_\theta$
has Lebesgue measure zero, and such that for all $\lambda >2$ and all
$\alpha\in S_\theta$, the Mathieu operator $H$ 
has an eigenvalue $E$ such that the corresponding eigenvector
$(u_n)$ is of exponential decay as $n\rightarrow \pm\infty$.}\par
\vskip.05in
\indent Moreover, she proved that 
the Mathieu operator has pure point spectrum for
$\lambda >2$, and further conjectured that the same conclusion
holds for all real $\alpha$. For comparison, for $\lambda =2$ there is purely singular
continuous spectrum; whereas for $0<\lambda <2$ there is purely
absolutely continuous spectrum.\par
\vskip.05in

\indent  We introduce the Hankel matrices 
$$\Gamma_c=[\cos\pi (\alpha +x\theta +k\theta )\,
u_{x+k}]_{x,k=1}^{\infty}, \quad 
\Gamma_s=[\sin\pi (\alpha +x\theta +k\theta )\, u_{x+k}]_{x,k=1}^\infty.\eqno(6.4)$$
\noindent There exists a bounded and measurable function $\Phi :{\bf
T}\rightarrow M_2({\bf C})$ such that 
$$\hat \Phi (k) =\left[ \matrix{\cos\pi (\alpha +k\theta )\,
u_k&\sin\pi (\alpha +k\theta )\,
u_k\cr \cos\pi (\alpha +k\theta )\,
u_k&\cos\pi (\alpha +k\theta )\,
u_k\cr}\right];\eqno(6.5)$$
\noindent the block Hankel operator associated with $\Phi$ is
$[\hat\Phi (j+k)]$, which becomes, after a rearrangement of the block
form, the matrix 
$$\Gamma_\Phi =\left[ \matrix{ \Gamma_c&\Gamma_s\cr
\Gamma_s&\Gamma_c\cr}\right].\eqno(6.6)$$
\noindent The negative Fourier coefficients of $\Phi$ are not uniquely
determined by $\Gamma_\Phi$, but may be chosen advantageously. Note that all of these operators are self-adjoint.\par 
\vskip.05in
\indent We introduce operators $K$ and $L$ by  
$$\left[\matrix{L&K\cr
K&L\cr}\right]=\left[\matrix{\Gamma_c^2+\Gamma_s^2
&\Gamma_c\Gamma_s+\Gamma_s\Gamma_c\cr
\Gamma_c\Gamma_s+\Gamma_s\Gamma_c&
\Gamma_c^2+\Gamma_s^2\cr}\right]=\Gamma_\Phi^2.\eqno(6.7)$$

\vskip.05in

\vskip.05in
\indent {\bf Theorem 6.2.} {\sl Let $u_n$ be as in Lemma 6.1. \par
\indent (i) The matrices of $K$ and $L$ are given by:
$$K(m,n)={{u_{m-1}u_{n}-u_{n-1}u_{m}}\over 
{2\lambda\sin \pi\theta
(m-n)}}\qquad (m\neq n; m,n\in{\bf N});
\eqno(6.8)$$
$$L(m,n)=\cos \pi\theta (m-n)\sum_{k=1}
^\infty u_{m+k}u_{n+k}\qquad
(m,n\in {\bf N}).$$  
\indent (ii) Then
$\Gamma_\Phi^2$ is a positive semidefinite trace-class operator. \par
\indent (iii) Further,
the eigenvectors of $K$, $\Gamma_c$ and $\Gamma_s$ are of exponential
decay as $x\rightarrow\infty$.\par
\indent (iv) There exists a
bounded and measurable function $\Psi_n :{\bf T}\rightarrow M_2({\bf C})$
such that the associated Hankel operator $[\hat \Psi_n (j+k)]$ has rank
less than or equal to $n$ and 
$$s_n(\Gamma_\Phi )=\Vert \Gamma_\Phi -\Gamma_{\Psi_n}\Vert =\Vert
\Phi-\Psi_n\Vert_\infty .\eqno(6.9)$$
\indent (v) The eigenvalues 
of $K,L, \Gamma_s$ and $\Gamma_c$ are of exponential decay.}\par
\vskip.05in
\indent {\bf Proof.} (i) First we observe that the formula for
$K(m,n)$ makes sense for $m\neq n$ since $\theta$ is irrational. 
The discrete Mathieu equation 
$$u_{n+1}+u_{n-1}+\lambda \cos 2\pi (n\theta +\alpha )\,
u_n=Eu_n\eqno(6.10)$$
\noindent gives the
system
$$\left[\matrix{u_n\cr u_{n+1}\cr}\right]=
\left[\matrix{0&1\cr -1 &E-\lambda\cos
2\pi (n\theta +\alpha )\cr}\right]\left[\matrix{u_{n-1}
\cr u_{n}\cr}\right].\eqno(6.11)$$
\noindent Writing $T(n)$ for the one-step transition matrix, we have
$${{T(n)^tJT(m)-J}\over{2\lambda \sin\pi\theta (m-n)}}=
-\left[\matrix{0&0\cr 0&1\cr}\right] \bigl( \sin\pi (n\theta +\alpha
)\cos \pi (m\theta +\alpha )+\cos \pi (n\theta +\alpha )\sin\pi
(m\theta +\alpha )\bigr);\eqno(6.12)$$
\noindent so with $a(m)=[u_{m-1}, u_{m}]^t$ we introduce
$$\tilde K(m,n)={{\langle Ja(m),
a(n)\rangle}\over{2\lambda\sin\pi\theta (m-n)}}\eqno(6.13)$$
\noindent which satisfies
$$\eqalignno{\tilde K(m+1,n+1)-\tilde K(m,n)&= \bigl( \sin\pi (n\theta +\alpha
)\cos \pi (m\theta +\alpha )\cr
&\quad +\cos \pi (n\theta +\alpha )\sin\pi
(m\theta +\alpha )\bigr)u_nu_m}$$
\noindent and $\tilde K(m+k,n+k)\rightarrow 0$ as $k\rightarrow\infty$. Hence by comparing entries of the products, we find
that $\tilde K=\Gamma_c\Gamma_s+\Gamma_s\Gamma_c.$\par
\indent (ii) This is clear, since the entries of the matrix of
$\Gamma_\Phi$ are real and summable.\par
\indent (iii) For any unit vector $\varphi$, we have 
$$\Gamma_c\varphi (x)=\sum_{k=0}^\infty \cos \pi (\alpha +\theta
x+\theta k)\, u_{x+k}\varphi_k\eqno(6.14)$$
\noindent so by the Cauchy--Schwarz inequality we have the uniform bound
$$\vert \Gamma_c\varphi (x)\vert\leq \Bigl(\sum_{k=x}^\infty 
u_k^2\Bigr)^{1/2},\eqno(6.15)$$
\noindent where the right-hand side decays exponentially as
$x\rightarrow\infty.$ A similar result applies with $\Gamma_s$, so
in particular the eigenvectors of $\Gamma_c$ and $\Gamma_s$ are of
exponential decay at infinity.\par
\indent Now let $(\varphi_j)_{j=1}^\infty$ be an orthonormal basis of
$\ell^2({\bf Z})$ consisting of eigenvectors of $\Gamma_s$ with
corresponding eigenvalues $\sigma_j$. Then  
$$\Gamma_c\Gamma_s\varphi (x)=\sum_{j=1}^\infty \sigma_j\langle \varphi_j,
\varphi \rangle \Gamma_c\varphi_j(x)\eqno(6.16)$$
\noindent where 
$$\sum_{j=1}^\infty \vert \sigma_j\langle \varphi_j,
\varphi \rangle\vert \leq \Bigl(\sum_{j=1}^\infty
\sigma_j^2\Bigr)^{1/2}\Bigl(\sum_{j=1}^\infty \langle \varphi
,\varphi_j\rangle^2\Bigr)^{1/2}=\Vert
\Gamma_s\Vert_{c^2}\Vert\varphi\Vert_{\ell^2}.\eqno(6.17)$$
\noindent This and a similar result for $\Gamma_s\Gamma_c$ imply $\vert
K\varphi (x)\vert$ decays exponentially as $x\rightarrow\infty.$\par
\indent (iv) This is immediate from the vectorial form of the matrical
AAK theorem [14].\par
\indent (v) Since $u_n$ decays exponentially as
$n\rightarrow\infty$, we can approximate the Hankel matrices with
finite matrices up to exponentially small error terms. For instance, we
can approximate $\Gamma_s$ by $\Gamma_s^{(N)}=[\sin\pi(\alpha +x\theta
+k\theta )u_{x+k}{\bf I}_{\{ (x,k): x+k\leq N\}}]$ which has rank less
than $N+1$ and the operator norm satisfies 
$$\Vert \Gamma_s-\Gamma_s^{(N)}\Vert\leq\sum_{k=N+1}^\infty k\vert
u_k\vert .\eqno(6.18)$$
\indent  The $s$-numbers satisfy
$$s_n(K)\leq s_n\Biggl(\left[\matrix{K&L\cr L&K\cr}\right]\Biggr)= 
s_n(\Gamma_\Phi^2)=s_n(\Gamma_\Phi )^2.\eqno(6.19)$$
\rightline{$\square$}\par
\vskip.05in
\indent A vectorial Hankel matrix $\Gamma_\Psi$ has finite rank if
and only if it has a rational symbol with coefficients of finite rank
by [14, p. 19].
Peller provides a formula for the rank in terms of the coefficients.\par
\vskip.05in
\indent {\bf Acknowledgements.} This research was partially supported by
EU Network Grant MRTN-CT-2004-511953 `Phenomena in High Dimensions'.
Andrew McCafferty's research was supported by EPSRC.\par
\vskip.1in
\noindent {\bf References}\par
\vskip.05in
\noindent 1. G. AUBRUN, A sharp small deviation inequality for the largest eigenvalue
of a random matrix, {\sl Springer Lecture Notes in Math.,} {\bf 1857}, (Springer, Berlin, 2005).\par 
\noindent 2. J. AVRON AND B. SIMON, Almost periodic Schr\"odinger
operators II. The integrated density of states, {\sl Duke Math. J.}
{\bf 50} (1983), 369--391.\par 
\noindent 3. G. BLOWER, {Operators associated with the soft and hard edges from unitary
ensembles}, {\sl J. Math. Anal. Appl.} {\bf 337} (2008), 239--265.
(doi:10.1016/j.jmaa.2007.03.084.)\par
\noindent 4. G. BLOWER, {Integrable operators and the squares of
Hankel operators}, {\sl J. Math. Anal. Appl.} {\bf 340} (2008),
943--953.\par
\noindent 5. A. BORODIN, A. OKOUNKOV AND G. OLSHANSKI, Asymptotics of
Plancherel measures for symmetric groups, {\sl J. Amer. Math. Soc.}
{\bf 13}
(2000), 481--515.\par
\noindent 6. J. BOURGAIN, On the spectrum of lattice Schr\"odinger
operators with deterministic potential. Dedicated to the memory of
Thomas H. Wolff, {\sl J. Anal. Math.} {\bf 87} (2002), 37--75.\par
\noindent 7. G.H. HARDY, J.E. LITTLEWOOD AND G. P\^OLYA,
{\sl Inequalities},
(Cambridge University Press, 1988)\par
\noindent 8. J.W. HELTON, Discrete time systems, operator models and
scattering theory, {\sl J. Funct. Anal.} {\bf 16} (1974), 15--38.\par
\noindent 9. R.A. HORN AND C.R. JOHNSON, {\sl Topics in Matrix
Analysis}, (Cambridge University Press, 1991).\par
\noindent 10. S. JITOMIRSKAYA, Metal-insulator transition for the almost
Mathieu operator, {\sl Ann. Math. (2)} {\bf 150} (1999), 1159--1175.\par
\noindent 11. K. JOHANSSON, Discrete orthogonal polynomial ensembles and
the Plancherel measure, {\sl Ann. Math. (2)} {\bf 153} (2001), 259--296.\par 

\noindent 12. W. MAGNUS AND S. WINKLER, {\sl Hill's equation}, (Dover
Publications, New York, 1966).\par
\noindent 13. A.N. MEGRETSKI\u I, V.V. PELLER AND S.R. TREIL, 
The inverse spectral problem for self-adjoint Hankel operators,
{\sl Acta Math.} 
{\bf 174} (1995), 241--309.\par
\noindent 14. V. PELLER, {\sl Hankel Operators and Their
Applications}, (Springer, New York, 2003).\par
\noindent 15. A.G. SOSHNIKOV, Determinantal random point fields,
2000,\par
\noindent {\sl arXiv.org:math/0002099.}\par
\noindent 16. C.A. TRACY AND H. WIDOM, {Fredholm determinants, differential
equations and matrix models}, {\sl Comm. Math. Phys.} {\bf 163} (1994), 33--72.\par
\noindent 17. E.T. WHITTAKER AND G.N. WATSON, {\sl A Course of Modern
Analysis}, fourth edition, (Cambridge University Press, 1965).\par 
\vfill
\eject
\end